\documentclass[preprint]{elsarticle}
\usepackage{lineno}

\usepackage{stmaryrd}
\usepackage{slashbox}
\usepackage{amsmath}
\usepackage{amssymb}
\usepackage{amsthm}
\usepackage{graphicx}
\usepackage{epic,eepic,epsfig}
\usepackage{color}
\usepackage{subfigure}  
\usepackage{placeins}   
\usepackage{multirow}
\usepackage{epstopdf}
\usepackage{varwidth}
\usepackage{float}
\usepackage[table]{xcolor}

\newtheorem{example}{Example}[section]

\newtheorem{remark}{Remark}[section]

\newtheorem{algorithm}{Algorithm}[section]

\usepackage{lineno} 

\newcommand{\BLbrace}{ \left\{\kern -0.36em\left\{     }
\newcommand{\BRbrace}{ \right\}\kern -0.36em\right\} }
\newcommand{\Lbrace}{ \left\{\kern -0.23em\left\{     }
\newcommand{\Rbrace}{ \right\}\kern -0.23em\right\} }

\newcommand{\BLbracket }{ \left[\kern -0.19em\left[    }
\newcommand{\BRbracket }{ \right] \kern-0.19em\right] }
\newcommand{\Lbracket }{ \left[\kern -0.09em\left[    }
\newcommand{\Rbracket }{ \right] \kern-0.09em\right] }

\allowdisplaybreaks[4]

\begin{document}

\begin{frontmatter}
	
\title{Iterative two-level algorithm for nonsymmetric or indefinite elliptic problems}

\author[SCNU]{Ming Tang}
\ead{mingtang@m.scnu.edu.cn}

\author[SCNU]{Xiaoqing Xing}
\ead{xingxq@scnu.edu.cn}

\author[GUET]{Ying Yang\corref{cor}}  
\ead{yangying@lsec.cc.ac.cn}

\author[SCNU]{Liuqiang Zhong}
\ead{zhong@scnu.edu.cn}

\cortext[cor]{Corresponding author}
\address[SCNU]{School of Mathematical Sciences, South China Normal University, Guangzhou 510631, China}

\address[GUET]{School of Mathematics and Computational Science, Guangxi Colleges and Universities Key Laboratory of Data Analysis and Computation,
	Guilin University of Electronic Technology, Guilin 541004, China}

\begin{abstract}
 In this paper, a new iterative two-level algorithm  is presented for solving the finite element discretization for nonsymmetric or indefinite elliptic problems. The iterative two-level algorithm uses the same coarse space as the traditional two-grid algorithm, but its ``fine space''   uses the higher oder finite element space under the coarse grid.
 Therefore, the iterative two-level algorithm only needs one grid, and the computational cost is much lower than the traditional iterative two-grid algorithm. Finally, compared with the traditional two-grid algorithm, numerical experiments show that the computational cost is lower to achieve the same convergence order.  
\end{abstract}


\end{frontmatter}


\section{Introduction}\setcounter{equation}{0}\label{Cha:1}  

We consider the following Dirichlet boundary value problems for general second-order elliptic partial differential equations
\begin{eqnarray}\label{PDEs:NSPD:2Order}
-\mathrm{div}(\boldsymbol{\alpha}(x) \nabla u)+\boldsymbol{\beta}(x) \cdot \nabla u+\gamma(x) u = f 
&\quad&
\mbox{in $\Omega$ }, 
\\ \label{DHBC:NSPD:2Order}
u = 0  &\quad& \mbox{on  $\partial \Omega$ },
\end{eqnarray}
where the coefficient $\boldsymbol{\alpha}(x) \in \mathbb{R}^{d} \times \mathbb{R}^{d}$ is smooth functions on $\bar{\Omega}$, and satisfies consistent ellipticity, namely, there are minimum and maximum eigenvalues $\alpha_{*}$ and $\alpha^{*}$ respectively, satisfying $\alpha_{*}|\xi|^2 \leqslant \mathbf{\alpha}(x)\xi \cdot \xi \leqslant \alpha^{*}|\xi|^2, ~~\forall \xi \in \mathbb{R}^d$. Both $\boldsymbol{\beta}(x)  \in \mathbb{R}^{d}$ and $\gamma(x) \in \mathbb{R}^1$ are smooth functions on $\bar{\Omega}$, and $f \in L^2(\Omega)$ is a given function. Noting that when $\boldsymbol{\beta}(x) \neq \boldsymbol{0}$, the continuous variational problems (CVP) of \eqref{PDEs:NSPD:2Order}-\eqref{DHBC:NSPD:2Order} are nonsymmetric, and when $\gamma(x)<0$, the CVP of \eqref{PDEs:NSPD:2Order}-\eqref{DHBC:NSPD:2Order} may be indefinite.

The two-grid (TG) algorithm was first introduced by Xu \cite{XuJC96:1759} for solving nonasymmetric or indefinite problems. 
The main idea of the TG algorithm  is to solve the original problems on the coarse mesh 
to obtain an approximate finite element solution, and then use the approximate solution to solve the corresponding linear symmetric positive definite problems on the fine mesh. 
Over the last two decades, the TG algorithm has been widely used to solve many problems, such as nonlinear elliptic problems \cite{BiCJWangC18:23,ZhongLQZhouLL21:105587}, semilinear parabolic equations \cite{MarionXu95:1170}, nonlinear
parabolic equations \cite{DawsonWheele94Meeting,DawsonWheele98:435}, Poisson-Nernst-Planck problems \cite{YangYLuBZ19:556}, and Maxwell equations \cite{ZhongLQLiuCM13:432,ZhongLQShuS13:93}.

The main objective of this paper is to propose a new iterative two-level algorithm. 
Compared with the traditional iterative two-grid algorithm, our algorithm needs only one mesh and takes less CPU time to achieve the same accuracy.

In this work, we first present elliptic problems discretized by the finite element method. Then, we propose an iterative two-level algorithm.  Finally, some numerical results are presented to illustrate the efficiency of the proposed algorithms.

\section{Finite element discretizations}\setcounter{equation}{0}
Given $S \subset \mathbb{R}^d(d=2, 3)$, we denote $W^{1, 2}(S)$ as the standard Sobolev space with norm $\| \cdot \|_{1, 2, S}$. For simplicity of notation,
we denote $\| \cdot \|_{1} = \| \cdot \|_{1, 2, \Omega}$, and $H^1_0(\Omega) := \{u\in H^1(\Omega) : u|_{\partial \Omega} =0 \}$ in the sense of trace. We define the following two bilinear forms
\begin{eqnarray} \label{Def:Bilinear:a}
a(u, v)&:=&\int_\Omega (\boldsymbol{\alpha}(x)\nabla u) \cdot \nabla v \mathrm{d}x, 
\\ 
\label{Def:Bilinear:hata}
\hat{a}(u, v)&:=&a(u, v)+\int_\Omega (\boldsymbol{\beta}(x) \cdot \nabla u+\gamma (x)u)v \mathrm{d}x.
\end{eqnarray}

Then, we obtain the CVP of  \eqref{PDEs:NSPD:2Order}-\eqref{DHBC:NSPD:2Order}: Find $u \in H_{0}^{1}(\Omega)$, such that
\begin{equation}\label{Eqn:2.10} 
\hat{a}(u, v)=(f, v),~~\forall v \in H_{0}^{1}(\Omega).
\end{equation}

We assume that $\Omega$ is partitioned by a quasi-uniform division $\mathcal{T}_h = \{\tau_i\}$.  By this we mean that $\tau_i$'s are simplexes of the size $h$ with $h \in (0,1)$ and $\bar{\Omega} = \cup_i \tau_i$. For the given quasi-uniform division $\mathcal{T}_h$, the conforming finite element space is defined as follows
\begin{equation*}\label{Eqn:2.1}
V_h^l:=\{ v \in C(\Omega) : v|_{\tau} \in \mathbb{P}_l(\tau),~ \forall \tau \in \mathcal{T}_h, ~v |_{\partial \Omega}=0 \},
\end{equation*}
where $\mathbb{P}_l(\tau)$ is the space of polynomial of degree not greater than a positive integer $l$ on the subdivision element $\tau$.
The discrete variational problems of \eqref{Eqn:2.10} is to find $u_h^l \in V_h^l$, such that
\begin{equation}\label{Eqn:2.11}
\hat{a}(u_h^l, v_h^l)=(f, v_h^l),~~\forall v_h^l\in V_h^l.
\end{equation}

\section{Iterative two-level  algorithm} \setcounter{equation}{0}
The basic mechanism of the classical iterative two-grid  algorithm is two quasi-uniform tetrahedral nested meshes of $\Omega$, namely the fine space $V_h^{k, l}$ and the coarse space $V_H^{k, l}$, with two different meshes sizes $h$ and $H(h<H)$. Furthermore, in the application given in the succeeding text, we shall always assume that
$
H=O\left(h^\lambda\right) \text { for some } 0<\lambda<1.
$

\begin{algorithm}[Algorithm 4.1 of  \cite{XuJC96:1759}]\label{Alg:Itg1}    
   Let $u_{h}^{l,0}=0$; assume that $u_{h}^{l,k}\in {V}_h^l(k\geq 0)$ has been obtained, $u_{h}^{l,k+1}\in {V}_h^l$ is defined as follows:
	\begin{description} 
		\item[1.]  Find $e_{H}^{l,k}\in V_H^l$, such that
		$
			\hat{a}(e_{H}^{l,k},v_H^l) = (f,v_H^l) - \hat{a}(u^{l,k}_{h},v_H^l), \ \forall v_H^l\in V_H^l.
		$
	
		\item[2.] Find $u_{h}^{l,k+1}\in V_h^l$, such that
		$
			a(u_{h}^{l,k+1},v_h^{l}) = 
			(f,v_h^l) - N( u_{h}^{l,k} + e_{H}^{l,k},v_h^l),\ \forall v^l_{h}\in V_h^l,
		$
		
		where $N(u, v):=\hat{a}(u, v)-a(u, v), \quad \forall u, v \in H_{0}^{1}(\Omega).$ 
	\end{description}
\end{algorithm}

Assume that $u_{h}^{l,k+1}\in V_h^l$ is the solution obtained by Algorithm \ref{Alg:Itg1} with $k\geq 1$, then we have (see Theorem 4.4 of \cite{XuJC96:1759})
\begin{equation}\label{Eroor}
	\| u - u_{h}^{l,k+1} \|_1\lesssim (h^l + H^{k+l}) \|u\|_{l+1},
\end{equation}
which means that the two-grid solution given by Algorithm \ref{Alg:Itg1}  can effectively approximate the finite element solution $u_h^l$ of   \eqref{Eqn:2.11}.

\begin{remark}
In order to obtain the optimal convergence order in \eqref{Eroor}, we should assume that  $H$  and   $h$  satisfy the relation $h=O(H^\frac{l+k}{l})$. 
However, there are exact nested grids, which satisfy $h=O(H^\frac{l+k}{l})$, very difficult to implement.
\end{remark}

Next, we consider replacing $V^l_h$ with $V_H^s(s\geq l+1)$ in the second step of Algorithm \ref{Alg:Itg1}, and obtain the following algorithm.

\begin{algorithm}\label{Alg:Itg2}
	   Let $\hat{u}_{H}^{l,0}=0$; assume that $\hat{u}_{H}^{s,k}\in {V}_H^s(k\geq 0)$ has been obtained, $\hat{u}_{H}^{s,k+1}\in {V}_H^s$ is defined as follows:
	\begin{description} 
		\item[1.]  Find $e_{H}^{l,k}\in V_H^l$, such that
		$
		\hat{a}(e_{H}^{l,k},v_H^l) = (f,v_H^l) - \hat{a}(\hat{u}^{s,k}_{H},v_H^l), \ \forall v_H^l\in V_H^l.
		$
		\item[2.] Find $\hat{u}_{H}^{s,k+1}\in V_H^s$, such that
		$
		a(\hat{u}_{H}^{s,k+1},v_H^{s}) = 
		(f,v_H^s) - N( \hat{u}_{H}^{s,k} + e_{H}^{l,k},v_H^s),\ \forall v^s_{H}\in V_H^s.
		$
	\end{description}
\end{algorithm}
\begin{remark}
Comparing Algorithm \ref {Alg:Itg1} with Algorithm \ref{Alg:Itg2}, 
although the SPD problems are solved in the second step of the algorithm, $V_H^{s}$ has fewer degrees of freedom than $V_h^{l}$ (See Table \ref{Tab:3}), which can greatly reduce the computation time.
\end{remark}

\begin{table}[!htbp]
	\renewcommand\arraystretch{1}
	\centering
	\begin{tabular}{cccccc}
		\hline
		$\mathcal{T}_H$  &  $dof(V_H^3)$&  $dof(V_h^3)$ &   $dof(V_H^4)$    &  $dof(V_H^5)$  & $dof(V_H^6)$\\ 
		\hline
		1/9  & 784   & 59536   & 1369  & 2116  & 3025 \\
		1/10 & 961 	 & 90601   & 1681  & 2601  & 3721 \\
		1/11 & 1156  & 132496  & 2025  & 3136 & 4489 \\
		1/12 & 1369  & 187489  & 2401  & 3721 & 5329\\
		\hline
	\end{tabular}
	\caption{ Taking $h=H^2$, $V_H^{4}$, $V_H^{5}$, $V_H^{6}$ and $V_h^{3}$ degrees of freedom comparison }\label{Tab:3}
\end{table}

\section{Numerical results}

In this section, numerical experiments are carried out to verify the effectiveness of the iterative two-level algorithm.  We performed all experiments for our iterative two-level algorithm with the help of the software package Fenics \cite{LoggMardal12Book}.

\begin{figure}[H]
		\centerline{\includegraphics[scale=1]{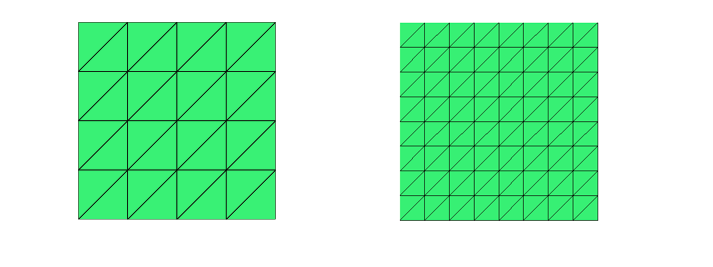}}
		\caption{Structured grids with $\mathcal{T}_{1/4}$ (left) and $\mathcal{T}_{1/8}$ (right).}
		\label{StructuredGrid:2D:4}
	\end{figure}

\begin{example}\label{example-1}
 
	We consider model problems \eqref{PDEs:NSPD:2Order}-\eqref{DHBC:NSPD:2Order}, where the computational domain is $\Omega=(0, 1)^2$,  the coefficients are $\alpha = 1.0$, $\boldsymbol{\beta}=(0, 0)^t$ and $\gamma = -10$, the exact solution  is $u= \sin(\pi x)\sin(\pi y)$, and $f$ can be obtained by substituting the exact solution into equation \eqref{PDEs:NSPD:2Order}.
\end{example}	
	We first partition the $x-$ axis and $y-$ axis of the domain  $\Omega$ into equally distributed $M$ subintervals, then divide each square into two triangles by using the line with slope $-1$. Hence, we obtain a sequence of nested and structured grids and  the corresponding meshes as $\mathcal{T}_{H}$ with $H=1/M$, where $M\geq 2$ is an integer, see Figure \ref{StructuredGrid:2D:4}. We choose the piecewise conform $l$ order finite element spaces $V_H^{l}$ based on the meshes $\mathcal{T}_{H}$. For Algorithm \ref{Alg:Itg1}, we choose $h=H^2$.		
 \begin{table}[!htbp]
	\centering
	\begin{tabular}{ccccccc}
		\hline
		\multicolumn{1}{l}{}  &\multicolumn{3}{c}{ Algorithm  \ref{Alg:Itg1}, $l=3$, $k=3$} &\multicolumn{3}{c}{Algorithm \ref{Alg:Itg2}, $l=3$, $s=6$, $k=3$}  \\ \cline{2-4} \cline{5-7}
		$H$  & $\|u-u_{h}^{3, 3} \|_1$ &$\|u-u_{h}^{3, 3} \|_1*H^{-6}$ & CPU &  $\|u-\hat{u}_{H}^{6, 3}\|_1$ & $\|u-\hat{u}_{H}^{6, 3}\|_1*H^{-6}$&CPU \\ \hline
		1/9  &9.8925E-07& 5.2573E-01& 6.1544 &5.7750E-08&	3.0691E-02&	0.2503\\
		1/10 &5.2609E-07& 5.2609E-01& 9.0338 &3.0706E-08&	3.0706E-02&	0.3082\\
		1/11 &2.9711E-07& 5.2635E-01& 16.177 &1.7339E-08&	3.0717E-02&	0.3644\\
		1/12 &1.7634E-07& 5.2655E-01& 24.047 &1.0290E-08&	3.0726E-02&	0.4347\\
		\hline
	\end{tabular}
	\caption{Compare the $H^1$ error estimate between the classical iterative two-grid  algorithm \ref{Alg:Itg1} and the iterative two-level algorithm \ref{Alg:Itg2}. }\label{lable-t3-7}
\end{table}

From Table \ref{lable-t3-7}, we can observe that both our algorithm and the traditional iterative mesh can reach the optimal convergence order. And to achieve the same accuracy, our algorithm uses less CPU time.

 \begin{table}[!htbp]
	\centering
	\begin{tabular}{ccccccc}
		\hline
		\multicolumn{1}{l}{}  &\multicolumn{3}{c}{ Algorithm  \ref{Alg:Itg2}, $l=3$, $s=4$, $k=3$} &\multicolumn{3}{c}{Algorithm \ref{Alg:Itg2}, $l=3$, $s=5$, $k=3$ }  \\ \cline{2-4} \cline{5-7}
		$H$  & $\|u-\hat{u}_{H}^{4, 3} \|_1$ &$\|u-\hat{u}_{H}^{4, 3} \|_1*H^{-4}$ & CPU &  $\|u-\hat{u}_{H}^{5, 3}\|_1$ & $\|u-\hat{u}_{H}^{5, 3}\|_1*H^{-5}$&CPU \\ \hline
		1/9  &3.6409E-05&	2.3888E-01&	0.1383&1.6093E-06&	9.5028E-02&	0.1825\\
		1/10 &2.3903E-05&	2.3903E-01&	0.1633&9.5141E-07&	9.5141E-02&	0.2172 \\
		1/11 &1.6334E-05&	2.3915E-01&	0.2785&5.9129E-07&	9.5228E-02&	0.2571 \\
		1/12 &1.1538E-05&	2.3924E-01&	0.2237&3.8298E-07&	9.5299E-02&	0.3020 \\
		\hline
	\end{tabular}
	\caption{ the $H^1$ error estimate of the iterative two-level  algorithm \ref{Alg:Itg2}. }\label{lable-t3-17}
\end{table}

From the Tables \ref{lable-t3-7}- \ref{lable-t3-17}, it can be found that when the degree of the coarse space polynomial degree $l=3$ is fixed, as the degree of the fine space polynomial increases once, the convergence order of the error estimates in $H^1$-norm for the solution of the algorithm \ref{Alg:Itg2} increase by one order. It can be seen that the errors in $H^1$-norm for the solution of the algorithm \ref{Alg:Itg2} depend on the value of $s$.
  
\begin{example}\label{example-2}
	
	We consider model problems \eqref{PDEs:NSPD:2Order}-\eqref{DHBC:NSPD:2Order}, where the computational domain is $\Omega=(0, 1)^2$,  the coefficients are $\alpha = 1.0$, $\boldsymbol{\beta}=(0, 0)^t$ and $\gamma = -10$, the exact solution  is $u= x(1-x)^2y(1-y)^2$, and $f$ can be obtained by substituting the exact solution into equation \eqref{PDEs:NSPD:2Order}.

\end{example}

\begin{table}[!htbp]
	\centering
	\begin{tabular}{ccccccc}
		\hline
		\multicolumn{1}{l}{}  &\multicolumn{3}{c}{ Algorithm  \ref{Alg:Itg1}, $l=3$, $k=3$} &\multicolumn{3}{c}{Algorithm \ref{Alg:Itg2}, $l=3$, $s=6$, $k=3$}  \\ \cline{2-4} \cline{5-7}
		$H$  & $\|u-u_{h}^{3, 3} \|_1$ &$\|u-u_{h}^{3, 3} \|_1*H^{-6}$ & CPU &  $\|u-\hat{u}_{H}^{6, 3}\|_1$ & $\|u-\hat{u}_{H}^{6, 3}\|_1*H^{-6}$&CPU \\ \hline
	     1/9  &6.0567E-08&	3.2188E-02&	6.0850&2.8919E-13&	1.5369E-07&	0.1901\\
		 1/10 &3.2255E-08&	3.2255E-02&	8.9463&1.2153E-13&	1.2153E-07&	0.2335 \\
		 1/11 &1.8235E-08&	3.2305E-02&	16.0683&7.9992E-14&	1.4171E-07&	0.2751  \\
		 1/12 &1.0832E-08&	3.2344E-02&	23.8708&7.8801E-14&	2.3530E-07&	0.3291 \\
		\hline
	\end{tabular}
	\caption{Compare the $H^1$ error estimate between the classical iterative two-grid  algorithm \ref{Alg:Itg1} and the iterative two-level algorithm \ref{Alg:Itg2}. }\label{lable-t3-117}
\end{table}

\begin{table}[!htbp]
	\centering
	\begin{tabular}{ccccccc}
		\hline
		\multicolumn{1}{l}{}  &\multicolumn{3}{c}{ Algorithm  \ref{Alg:Itg2}, $l=3$, $s=4$, $k=3$} &\multicolumn{3}{c}{Algorithm \ref{Alg:Itg2}, $l=3$, $s=5$, $k=3$ }  \\ \cline{2-4} \cline{5-7}
		$H$  & $\|u-\hat{u}_{H}^{4, 3} \|_1$ &$\|u-\hat{u}_{H}^{4, 3} \|_1*H^{-4}$ & CPU &  $\|u-\hat{u}_{H}^{5, 3}\|_1$ & $\|u-\hat{u}_{H}^{5, 3}\|_1*H^{-5}$&CPU \\ \hline
		 1/9  &1.9981E-06&	1.3110E-02&	0.1035 &5.2140E-08&	3.0788E-03&	0.1380\\
		 1/10 &1.3129E-06&	1.3129E-02&	0.1214 &3.0796E-08&	3.0796E-03&	0.1637\\
		 1/11 &8.9783E-07&	1.3145E-02&	0.1420 &1.9126E-08&	3.0803E-03&	0.1938\\
		 1/12 &6.3456E-07&	1.3158E-02&	0.1661 &1.2381E-08&	3.0809E-03&	0.2273\\
		\hline
	\end{tabular}
	\caption{ the $H^1$ error estimate of the iterative two-level algorithm \ref{Alg:Itg2}. }\label{lable-t3-177}
\end{table}

The same conclusions can be observed in Tables \ref{lable-t3-117} and \ref{lable-t3-177} as in Tables \ref{lable-t3-7} and \ref{lable-t3-17}.

\section*{Acknowledgment}

The first, second and fourth authors are supported by the National Natural Science Foundation of China (No. 12071160).  The  second  author is also supported by the National Natural Science Foundation of China (No. 11901212). The third author is supported by the National Natural Science Foundation of China (No. 12161026), Guangxi Natural Science Foundation (No. 2020GXNSFAA159098).

\section*{References}

\bibliographystyle{abbrv}


\end{document}